\documentclass[11pt, reqno]{amsart}
\usepackage{lmodern}
\usepackage{amsmath, amsthm, amssymb, amsfonts}
\usepackage[normalem]{ulem}
\usepackage{hyperref}
\usepackage[all,cmtip]{xy}
\usepackage{verbatim}
\usepackage{nccmath}
\usepackage{stmaryrd}
\usepackage{caption}
\usepackage{mathrsfs}
\usepackage{mathtools}
\usepackage{esvect}
\usepackage{cite}
\usepackage{bbm}
\usepackage{eucal}

\usepackage{mathbbol}
\usepackage{tabularx}
\usepackage[toc,page]{appendix}

\usepackage{tikz-cd}

\theoremstyle{plain}
\newtheorem{thm}{Theorem}[section]
\newtheorem{cor}[thm]{Corollary}

\newtheorem{lemma}[thm]{Lemma}

\newtheorem*{conjno}{Conjecture}

\renewcommand{\arraystretch}{2}

\theoremstyle{definition}
\newtheorem{defn}[thm]{Definition}
\newtheorem{example}[thm]{Example}
\makeatletter
\newcommand\ackname{Acknowledgements}
\if@titlepage
\newenvironment{acknowledgements}{%
	\titlepage
	\null\vfil
	\@beginparpenalty\@lowpenalty
	\begin{center}%
		\bfseries \ackname
		\@endparpenalty\@M
\end{center}}%
{\par\vfil\null\endtitlepage}
\else

\fi
\makeatother

\theoremstyle{remark}

\newcommand{\BC}{{\mathbb{C}}}

\newcommand{\BF}{{\mathbb{F}}}

\newcommand{\BQ}{{\mathbb{Q}}}
\newcommand{\BR}{{\mathbb{R}}}

\newcommand{\BT}{{\mathbb{T}}}

\newcommand{\BZ}{{\mathbb{Z}}}

\newcommand{\CC}{{\mathcal C}}

\newcommand{\CF}{{\mathcal F}}

\newcommand{\CI}{{\mathcal I}}

\newcommand{\CL}{{\mathcal L}}
\newcommand{\CM}{{\mathcal M}}
\newcommand{\CN}{{\mathcal N}}
\newcommand{\CO}{{\mathcal O}}

\newcommand{\CS}{{\mathcal S}}
\newcommand{\CT}{{\mathcal T}}

\newcommand{\FM}{{\mathfrak{M}}}

\newcommand{\vsf}{{\mathsf{v}}}
\newcommand{\wsf}{{\mathsf{w}}}

\newcommand{\rsf}{{\mathsf{r}}}
\newcommand{\dsf}{{\mathsf{d}}}

\newcommand{\asf}{{\mathsf{a}}}
\newcommand{\msf}{{\mathsf{m}}}

\newcommand{\Usf}{{\mathsf{U}}}

\newcommand{\ch}{{\mathrm{ch}}}
\DeclareMathOperator{\Hilb}{Hilb}
\newcommand{\rk}{{\mathrm{rk}}}

\newcommand{\Cohc}{\mathrm{Coh}}

\newcommand{\Ms}{\check{M}}
\newcommand{\Mp}{\hat{M}}
\newcommand{\Msp}{\hat{\check{M}}}

\newcommand{\QMs}{\check{\mathsf{QM}}}
\newcommand{\QMp}{\hat{\mathsf{QM}}}
\newcommand{\GWs}{\check{\mathsf{GW}}}

\newcommand{\QMsp}{\hat{\check{\mathsf{QM}}}}

\newcommand{\Nsp}{\hat{\check{N}}}
\newcommand{\Thetas}{\check{\Theta}}
\newcommand{\Thetap}{\hat{\Theta}}
\newcommand{\Thetasp}{\hat{\check{\Theta}}}
\newcommand{\alphas}{\check{\alpha}}
\newcommand{\alphap}{\hat{\alpha}}
\newcommand{\alphasp}{\hat{\check{\alpha}}}

\newcommand{\tr}{{\mathrm{tr}}}

\DeclareFontFamily{OT1}{rsfs}{}
\DeclareFontShape{OT1}{rsfs}{n}{it}{<-> rsfs10}{}
\DeclareMathAlphabet{\curly}{OT1}{rsfs}{n}{it}

\newcommand\Hom{\operatorname{Hom}}
\newcommand{\p}{\mathbb{P}}

\newcommand{\Mbar}{{\overline M}}

\newcommand{\td}{\mathrm{td}}

\newcommand{\Pic}{\mathop{\rm Pic}\nolimits}

\newcommand{\ev}{{\mathrm{ev}}}
\newcommand{\evsf}{{\mathsf{ev}}}

\newcommand{\br}{\mathsf{br}}
\newcommand{\Aut}{\mathrm{Aut}}

\newcommand{\FMs}{\check{\mathfrak{M}}}

\newcommand{\FMp}{\hat{\mathfrak{M}}}
\newcommand{\FMsp}{\hat{\check{\mathfrak{M}}}}

\begin{document}
	
	\title[Quasimaps to moduli spaces of sheaves and related topics]
	{Quasimaps to moduli spaces of sheaves and related topics}
	
	\author{Denis Nesterov}
	\address{University of Vienna, Faculty of Mathematics}
	\email{denis.nesterov@univie.ac.at}
	\maketitle
	\begin{abstract} This is a review article based on a mini-course composed of four talks given by the author at UC Berkeley.

	\end{abstract}
	\setcounter{tocdepth}{1}
\tableofcontents\begin{center}
	
\end{center}
\section{Introduction}
 The purpose of this article is to explain various facets of a relation between 
 \begin{itemize}
 	\item  Donaldson--Thomas theory of a threefold $S\times C$,
 	\item Vertex function (also known as $I$-function),
 	\item Gromov--Theory of a moduli space of sheaves $M$ on the surface $S$.
 	\end{itemize}
 The relation is summarized in the figure below.  Quasimaps play an essential role in uncovering this relation, since it can be seen as a quasimap wall-crossing. We also discuss related phenomena in the order indicated in the table of contents. Applications can be found in Section \ref{applications1}, \ref{applications2} and \ref{applications3}.
	\begin{figure} [h!]	\vspace{1.5cm}
	\scriptsize
	\[
	\begin{picture}(200,0)(-30,-30)
		\thicklines
		\put(24,0){\line(1,0){31}}
		\put(95,0){\line(1,0){31}}
		\put(75,0){\makebox(0,0){\textsf{DT Vertex}}}
		\put(143,0){\makebox(0,0){ $\mathsf{GW}_C(M)$}}
		\put(0,0){\makebox(0,0){$\mathsf{DT_{rel}}(S\times C)$}}
	\end{picture}
\vspace{-0.5cm}	\]
	\caption{Summary}
	\label{summary}
\end{figure}

\noindent \textbf{Acknowledgments.} I am infinitely grateful to Andrei Okounkov and Mina Aganagi\'c for inviting me to Berkeley and giving me opportunity to speak. I also thank everyone  I met in Berkeley for such a wonderful month. 
\section{Quasimaps to moduli spaces of sheaves}

\subsection{\underline{Gro}mov and \underline{Gro}thendieck}
Let $S$ be a projective surface and $C$ be a smooth curve over $\BC$.  Assume for simplicity that $b_1(S)=0$.  A moduli space of non-constant maps 
\[ M_C(S)=\{f \colon C \rightarrow S\}\]
is not compact in general. For example, it fails to be compact already  for maps of degree 2 from $\p^1$ to $\p^2$. There exist various compactifications of this moduli space, two of which will be of interest for us. 
\\

\begin{figure}[h!]
	\centering
	\includegraphics[scale=0.4]{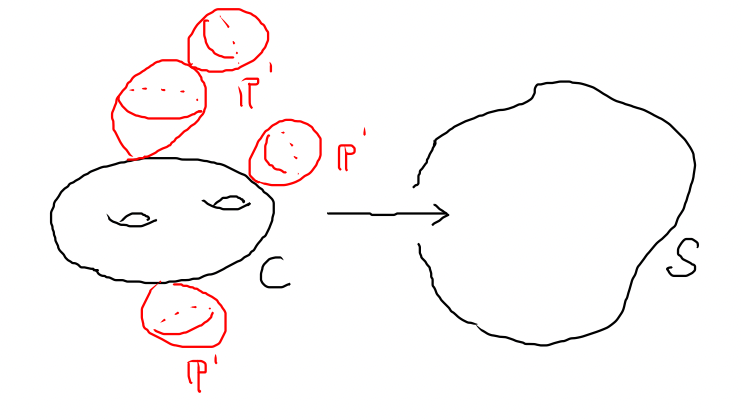}
	\caption{Maps from bubbling $C$}
	\label{bubbles}
\end{figure}

\textbf{Stable maps.} The first compactification is provided by maps from \textit{bubbling} $C$.  Bubbling $C$ is a curve $\tilde C$ with trees of $\p^1$ attached to $C$ at nodes,  see Figure \ref{bubbles}. We then consider \textit{stable maps} from these bubbling curves,   
 \[ \Mbar_C(S)=\left\{ f \colon \tilde C \rightarrow S \mid \deg(f_{|\mathrm{\p^1}})>0 \right\} /\hspace{-0.1cm}\sim,\] 
such that maps are identified by automorphisms of bubbles.  In the symplectic category,  compactness of $\Mbar_C(S)$  is a consequence of Gromov compactness theorem \cite{Gro}.  In the complex algebraic category, this compactification was constructed by Kontsevich \cite{Kon}. 

\begin{figure}[h!] \vspace{0.3cm}
	\centering
	\includegraphics[scale=0.4]{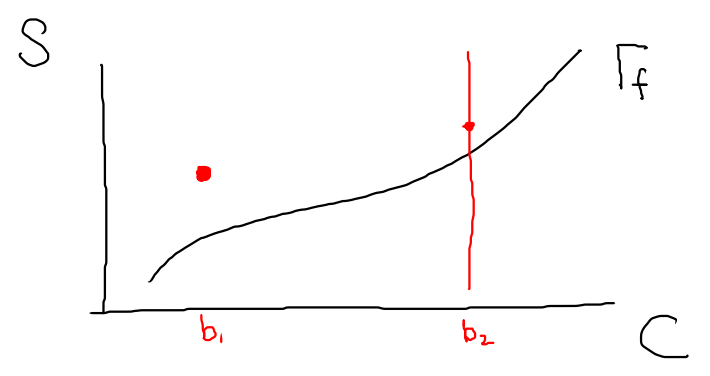}
	\caption{Degenerate graphs}
	\label{Degenerate}
	\vspace{-0.3cm}
\end{figure}

\vspace{-0.3cm}
\textbf{Degenerate graphs.} The second compacification is provided by Hilbert schemes constructed by Grothendieck \cite{Groth}. By associating to a map $f$ its graph $\Gamma_f \subset S\times C$, we obtain an embedding  

\begin{align*}
M_C(S) &\hookrightarrow \Hilb_1(S\times C), \quad f \mapsto \Gamma_f,
\end{align*}
where $\Hilb_1(S\times C)$ is a Hilbert scheme of one-dimensional subschemes on $S\times C$ of degree 1 over $C$ (the subscript indicates the degree, not the dimension). Hilbert schemes will contain one-dimensional subschemes which are not graphs of maps due to the presence of floating points and vertical components, see Figure \ref{Degenerate}. Compactness of this moduli space follows from the general theory of Hilbert schemes and relies on projectivity. \pagebreak

It might seem that these two compactifications are very different in nature, but this is not true - we just need to change the point of view. 
\subsection{Quasimaps to $S$}\label{points} Let us define another space 
\[\CS=\{\text{sheaves $F$ on $S$ }\mid  \ch(F)=(1,0,-1) \}.\]
Points of $\CS$ consist of the following three types of sheaves:
\begin{itemize}
	\item[(1)] ideal sheaves of points $I_p$;
	\item[(2)]  extensions of  ideal sheaves of points  by structure sheaves of points, e.g.\ $I_{p_1 \cup p_2} \oplus \CO_{p_2};$
	\item[(3)] extensions of ideal sheaves  of curves by structure sheaves of curves, e.g.\ $I_{\Gamma \cup p} \oplus \CO_{\Gamma}.$ 
	\end{itemize}
In particular, by existence of (1), there is a natural embedding $S \hookrightarrow \CS$. 
By existence of (2) and (3), one can readily see that $\CS$ is immensely bigger than $S$, it is not even bounded. Its usefulness in this context is due to quasimaps, which  were introduced by Ciocan-Fontanine--Kim--Maulik for certain GIT stacks \cite{CFKM}.
\begin{defn} A \textit{quasimap} to $S$ is a map $f \colon C \rightarrow \CS$ mapping generically to $S\subset \CS$. We denote quasimaps by $f\colon C \dasharrow  S$.
\end{defn}

\begin{lemma} \label{graph} Treating $\Hilb_1(S\times C)$ as a moduli space of ideal sheaves, there exists a natural identification, 
	\begin{align*}
	\Hilb_1(S\times C) &\xrightarrow{\sim } \{ f\colon C \dasharrow S \}, 
	\end{align*}
such that the quasimap $f$ associated to an ideal $I$ is defined by $f(p):=I_{|S\times p}$. 
	\end{lemma}
\textit{Proof.} This is almost a tautology. By the construction of $\CS$, we have
\begin{align*}
	f\colon C\rightarrow \CS &\iff \CF\text{ on } S\times C \text{ flat over \textit{C} with} \det(\CF)=\CO,\\
	f \text{ maps generically to } S &\iff \CF \text{ is torsion-free.} 
	\end{align*}
A rank 1 torsion-free sheaf with $\det(\CF)=\CO$ on a 3-fold is an ideal sheaf of a subscheme. We invite the reader to fill out the rest of the details. \qed 
\\


Given a quasimap $f \colon C \dasharrow S$, where do points that do not map to $S$ go? Applying the identification from above and looking at Figure \ref{Degenerate}, we see that 
\begin{align*}
& b_1 \mapsto \text{points of type (2)} \\
 & b_2 \mapsto \text{points of type (3)}.
 \end{align*}

Analogously, we can define
\[\CS^{[d]}=\{\text{sheaves $F$ on $S$ }\mid  \ch(F)=(1,0,-d) \}.\]
There is an embedding of Hilbert scheme of $d$-points,  $S^{[d]} \hookrightarrow \CS^{[d]}.$ 
By exactly the same arguments, we obtain a natural identification 
\[\Hilb_d(S\times C) \cong \{ f\colon C \dasharrow S^{[d]}  \},\]
where as before the subscript indicates that we consider subschemes of degree $d$ over $C$. 

\subsection{Quasimaps to moduli spaces of sheaves} 
Consider the lattice of algebraic classes on $S$ with the intersection pairing,
\[\Lambda := H^{\mathrm{alg}}(S), \quad  (a, b)= \int_S a\cdot b.\] 
 We fix a Chern character $\vsf \in \Lambda$ and an ample line bundle $\CO_S(1) \in \mathrm{Amp}(S)$. Define the associated moduli spaces 
 \begin{align*}
 \FM&(\vsf)=\{\text{sheaves } F \text{ on } S \mid \ch(F)=\vsf \}  \\
&\rotatebox[origin=c]{90}{$\subset$}   \\ 
 M&(\vsf)=\{ \text{Gieseker-semistable sheaves for } \CO_S(1)\}.
 \end{align*}
 We require the following:
\begin{itemize}
	\item semistable\ =\ stable,
	\item there exists a universal family\footnote{ In practice, one should require existence of $G \in K_0(S)$ such that $(\vsf,\ch(G)\cdot \td_S)=1$. } $\BF$ on $\FM(\vsf)$. 
	\end{itemize}
The second assumption can be dropped at expense of working with a finite gerbe; the first assumption can be also dropped in certain situations, see Section \ref{d0}.

Once everything is carefully defined, we obtain a copy of $\Lambda_\vsf:=\vsf^\perp$ sitting inside the cohomology of $\FM(\vsf)$,
\[\Lambda_\vsf \subset H^2(\FM(\vsf)),\]
we will thereby  treat $\Lambda_\vsf$ as divisor classes on $\FM(\vsf)$.  A choice of the universal family $\BF$ gives a section of the inclusion $\Lambda_\vsf \hookrightarrow \Lambda $,  which therefore gives a dual inclusion $\Lambda_\vsf^\vee \hookrightarrow \Lambda^\vee$. 
\begin{defn}
	Let $(C,p_i)$ be a marked nodal curve. A \textit{quasimap} 
	\[f \colon C \dashrightarrow M(\vsf)\]
	 is a map $f \colon C \rightarrow \FM(\vsf)$, such that 
	 \begin{itemize}
	 	\item $f$ maps generically to $M(\vsf)$,
	 	\item $f(\text{nodes}, p_i) \subset M(\vsf)$.
	 	\end{itemize}
Using the embedding $\Lambda_\vsf \subset H^2_\lambda(\FM(\vsf))$, we say $f$ is of \textit{degree} $\beta \in \Lambda^\vee_\vsf \subset \Lambda^\vee$, if 
 	\begin{itemize}
 		\item $\deg(f^*\CL) =\beta(\CL).$ 
 		\end{itemize}
  We call the finite set of points $\{p \in C \mid f(p) \in \FM(\vsf)\setminus M(\vsf)\}$ \textit{base points}, see Figure \ref{Degenerate}.
 	\end{defn}
 Stability of quasimaps slightly differs from stability of maps, as we have base points which we should exploit. In fact, the base points allow us to get rid of \textit{rational tails}. Rational tails of a curve are $\p^1$-components with one special point\footnote{Either nodes or marked points.}, also known as \textit{bubbles}, see Figure \ref{bubbles}. 
 \begin{defn} \label{stableq} A quasimap $f$ is \textit{stable}, if 
 	\begin{itemize}
 		\item $C$ does not have rational tails,
 		\item $|\Aut(f)| <\infty$.
 		\end{itemize}
 	\end{defn}
 We define a moduli space of stable quasimaps,
 \[Q_{g,n}(M(\vsf),\beta)=\left\{\text{stable }f\colon C \dasharrow M(\vsf) \ \Bigl\rvert \ \arraycolsep=0.1pt\def\arraystretch{1} \begin{array}{c} \deg(f)=\beta, \\[.001cm] g(C)=g, |p_i|=n \end{array} \right\}.\]
 Since $\FM(\vsf)$ is not bounded, the following theorem is a bit unexpected. 
 \begin{thm}[ \hspace{-0.1cm}\cite{N}] \label{main} $Q_{g,n}(M(\vsf),\beta)$ is a proper Deligne-Mumford stack. If $M(\vsf)$ is a smooth, then $Q_{g,n}(M(\vsf),\beta)$ is quasi-smooth. 
 	\end{thm}
 \textit{Sketch of Proof.} The proof boils down to the following three independent properties. 
 \begin{itemize}
 	\item  \textit{Positivity.} There exists a line bundle $\CL$ on $\FM(\vsf)$, such that for all non-constant quasimaps we have
 	\[\deg(f^*\CL)>0.\]
 	\item \textit{Hartogs' property.} Given a family of projective nodal curves $\CC \rightarrow \Delta$ over a discrete valuation ring $\Delta$  and a quasimap defined outside of a regular closed point $p\in \CC$, 
 	\[ f^\circ \colon \CC \setminus p \colon f  \rightarrow \FM(\vsf),\]
 	 then $f^\circ$  extends uniquely,  $f \colon \CC  \rightarrow \FM(\vsf)$. 
 	 \item \textit{Quasi-smoothness.} The stack $Q_{g,n}(M(\vsf),\beta)$ has a natural derived enchacement provided by Lurie's derived mapping stacks with the target $\BR\FM(\vsf)$, see \cite{Lur}. By the argument from \cite[Proposition 3.12]{Nqm}, the derived enhancement is quasi-smooth, if $M(\vsf)$ is smooth.  The virtual tangent complex is given as follows
 	 \[ \BT^{\mathrm{vir}}_{Q_{g,n}(M(\vsf),\beta)}= \pi_*f^*\BT^{\mathrm{vir}}_{\FM(\vsf)}. \]
 	\end{itemize}
 \qed
\subsection{Quasimaps and sheaves} 
 Let $\CC_{g,n} \rightarrow \CM_{g,n}$ be the universal curve over the moduli space of semistable curves (curves without rational tails). By evoking the argument of Lemma \ref{graph}, we can identify a moduli space of quasimaps with a moduli space of sheaves on $S\times C$ with varying $C$. 
 \begin{thm}[\vspace{-0.2cm}\cite{N}]\label{ident} There is a natural identification, 
 	\begin{align*}
 	Q_{g,n}(M(\vsf),\beta) \xrightarrow{\sim} M_\vsf(S\times \CC_{g,n}/ \CM_{g,n}, \beta^\vee), \quad f \mapsto \CF,
 	\end{align*}
 such that $\beta^\vee $ is the class dual to $\beta \in \Lambda^\vee$ with respect to the intersection pairing and 
 \[ \ch(\CF)= (\vsf, \beta^\vee ) \in \Lambda  \oplus \Lambda(-2)= H^{\mathrm{alg}}(S\times C).\] 
 	\end{thm}
 This identification can be taken as a dentition of the space on the right, examples below will illustrate its most important features.  

\begin{example} If $M(\vsf)=S^{[d]}$, then 
	\[ 	Q_{g,n}(S^{[d]},\beta)  \cong \Hilb_d(S\times \CC_{g,n}/ \CM_{g,n}, \beta^\vee).\]
If we specialise to a fixed smooth curve $C$, then 
		\[ Q_{C}(S^{[d]},\beta) \cong \Hilb_d(S\times C, \beta^\vee).\] 
If we choose a marked curve $(C,p)$, then we obtain a moduli space of ideals \textit{relative} to a divisor $S_p:=S\times p \subset S\times C$,
		\[ Q_{(C,p)}(S^{[d]} \cong \Hilb_d(S\times C/ S_p, \beta^\vee).\]
Relative stability of ideal sheaves  becomes equivalent to stability of quasimaps with a marking. 
	\end{example}
\begin{example} \label{exmp1}
More generally, for an arbitrary $M(\vsf)$, we have 
 
	\[Q_{C}(M(\vsf),\beta)  \cong M_
\vsf(S\times C,\beta^\vee),\]
	where $M_\vsf(S\times C,\beta^\vee)$ is the moduli space of sheaves stable with respect to $\CO_S(1)\boxtimes\CO_C(k)$ for $k\gg 0$. The determinant line bundle is fixed but in an exotic way, we refer to \cite[Lemma 3.15]{N} for more details on that. For a marked curve $(C,p)$, we obtain 
	\[ Q_{(C,p)}(M(\vsf),\beta) \cong M_\vsf(S\times C/ S_p, \beta^\vee),\]
where  $M_\vsf(S\times C/ S_p, \beta^\vee)$ is a moduli space of sheaves stable with respect to $\CO_S(1)\boxtimes\CO_C(k)$ for $k\gg 0$, whose restriction is also stable on $S_p$.  
	\end{example}

\subsection{Perverse quasimaps}
By changing the heart of $D^b(S)$, we get  quasimaps of different flavour. For example, consider the torsion pair
\begin{align*}
\CT&=\{A \in \Cohc(S) \mid \dim(A)=0\}  \\
\CT^{\perp}&=\{B \in \Cohc(S) \mid \Hom(A,B)=0 \ \forall A \in \CT\}.
\end{align*}
We then define the tilted heart associated to the torsion pair, 
\[\Cohc(S)_{\#}:=\langle \CT^\perp, \CT[-1] \rangle.\]
Let
\[ \CS^{[d]}_{\#}=\{ \text{objects $F$ in  }\Cohc(S)_{\#} \mid  \ch(F)=(1,0,-d) \}, \] 
be a moduli space of objects on $\Cohc(S)_{\#}$. The Hilbert scheme of $d$-points sits inside, 
\[ S^{[d]} \subset \CS^{[d]}_{\#}.\]
This construction puts a leash on the points of type (2). The moduli space of quasimaps to the pair $S^{[d]} \subset \CS^{[d]}_{\#}$, denoted by $Q_{g,n}(S^{[d]},\beta)_{\#}$,  is isomorphic to the moduli spaces of stable pairs of \cite{PT}, 
\[Q_{g,n}(S^{[d]},\beta)_{\#} \cong P_d(S\times \CC_{g,n}/ \CM_{g,n},\beta^\vee).\]

\section{Wall-crossing, Vertex and applications}
We want to compare stable quasimaps with stable maps. To do so, we will introduce a stability that interpolates between the two, allowing both base points and rational tails up to a given degree. 
\subsection{$\epsilon$-stability}
\begin{defn}
	Given a point $b \in C$ of a quasimap $f\colon C \dasharrow M(\vsf)$. We define another quasimap,
	\[f_b \colon C \dasharrow M(\vsf),\]
	with the following properties:
	\begin{itemize}
		\item $f_b=f$ on $C \setminus b$, 
		\item $f(b) \in M(\vsf)$. 
		\end{itemize}
The quasimaps $f_b$ is called \textit{stabilisation} of $f$ at $b$. It exists by properness\footnote{In fact, one just needs Langton's semistable reduction, which holds in greater generality, e.g. for sheaves on quasi-projective surfaces.} of $M(\vsf)$. For the following theorem, we use the line bundle that is briefly mentioned in Theorem \ref{main}. 

	\end{defn}
\begin{thm} There exists $\CL \in \Pic(\FM(\vsf))$, such that 
	\begin{itemize}
	\item $\forall$ non-constant quasimaps, $\deg_\CL(f):= \deg(f^*\CL)>0$,
	\item $\forall$ base points $b\in C$, $\deg_\CL(b):=\deg(f^*\CL)-\deg(f_b^*\CL)>0$. 
	\end{itemize}
	\end{thm}
\begin{defn} Given $\epsilon \in \BR_{>0}$. A quasimap $f$ is $\epsilon$-stable, if 
	\begin{itemize}
		\item $\forall$ rational tails $R, \  \deg_\CL(f^*\CL_{|R})>1/\epsilon$,
		\item $\forall$ base points $b, \ \deg_\CL(b)\leq 1/\epsilon$,
		\item $|\Aut(f)|<\infty$. 
		\end{itemize}
	\end{defn}
We introduce the following notation, 
\begin{align*}
\epsilon&=+ \quad \text{ if }\epsilon\gg 1 \\
\epsilon&=- \quad \text{ if } \epsilon\ll 1. 
\end{align*}
One can readily verify that $\epsilon$-stability specialises to stable maps and stable quasimaps, 
\begin{align*}
	+\text{-stable quasimaps } &= \text{ stable maps } \\
	- \text{-stable quasimaps } &= \text{ stable quasimaps }
	\end{align*}
We define moduli spaces of $\epsilon$-stable quasimaps, 
\[Q^{\epsilon}_{g,n}(M(\vsf),\beta)=\left\{\epsilon\text{-stable }f\colon C \dasharrow M(\vsf) \  \Bigl\rvert \ \arraycolsep=0.1pt\def\arraystretch{1} \begin{array}{c} \deg(f)=\beta, \\[.001cm] g(C)=g, |p_i|=n \end{array} \right\}.\]
\begin{thm} A moduli space $Q^{\epsilon}_{g,n}(M(\vsf),\beta)$ is a proper Deligne-Mumford stack. If $M(\vsf)$ is smooth, then $Q^{\epsilon}_{g,n}(M(\vsf),\beta)$ is quasi-smooth. 
	\end{thm}
\textit{Proof.} See Theorem \ref{main}. \qed

\subsection{Quasimap invariants}
There exist usual structures needed to define Gromov-Witten type invariants: 
\begin{itemize}
	\item virtual fundamental class
	\[[Q^{\epsilon}_{g,n}(M(\vsf),\beta)]^{\mathrm{vir}}\in H_{\mathrm{vdim}}(Q^{\epsilon}_{g,n}(M(\vsf),\beta)); \]
	\item evaluation maps 
	\[\ev_i: Q^{\epsilon}_{g,n}(M(\vsf),\beta) \rightarrow M(\vsf), \quad f \mapsto f(p_i);\]
\item $\psi$-classes, 
\[ L_{i}=T^*_{p_i}C , \quad \psi_i=\mathrm{c}_1(L_i). \]
	\end{itemize}
By Example \ref{exmp1}, with respect to the identification
	\[ Q_{g,n}(M(\vsf),\beta) \xrightarrow{\sim} M_\vsf(S\times \CC_{g,n}/\CM_{g,n},\beta^\vee),\]
 markings correspond to relative divisors, while evaluation maps correspond to evaluations at these relative divisors $\CF \mapsto \CF_{|S_{p_i}}$. 
\begin{defn} For classes $\lambda_i \in H^*(S^{[d]})$, we define $\epsilon$-stable quasimap invariants 
	\[\langle \lambda_1 \psi^{k_1}_1, \dots,\lambda_n \psi^{k_n}_n \rangle^{\epsilon}_{g,\beta}= \int_{[Q^{\epsilon}_{g,n}(M(\vsf),\beta)]^{\mathrm{vir}}} \prod^{n}_{i=1} \ev^*_i \lambda_i \cdot \psi^{k_i}.\]
	\end{defn}
\begin{example}
	 $\langle \lambda_1 \psi^{k_1}_1, \dots,\lambda_n \psi^{k_n}_n \rangle^{+}_{g,\beta}$ are Gromov-Witten invariants of $M(\vsf)$. 
	\end{example}
\begin{example} $\langle \lambda_1 \psi^{k_1}_1, \dots,\lambda_n \psi^{k_n}_n \rangle^{-}_{g,\beta}$ are Donaldson-Thomas invariants of $S\times \CC_{g,n}/ \CM_{g,n}$ with relative insertions. More specifically, if $(g,n)=(0,3),$ then $ \langle \lambda_1,  \lambda_2 ,\lambda_3  \rangle^{-}_{0,\beta}$ is a Donaldson-Thomas invariant of $S\times \p^1$ relative to $S_{0,1,\infty}$ with relative insertions. 
	\end{example}
\subsection{Vertex}  We define the following space, 
\[ V(M(\vsf), \beta) := \{ f \colon \p^1 \dasharrow M(\vsf) \mid \deg(f)=\beta, \ f(\infty) \in M(\vsf)\},\] 
which by definition has an evaluation map 
\[ \evsf \colon  V(M(\vsf),\beta) \rightarrow M(\vsf), \ f \mapsto f(\infty).\]
 Let the torus $\BC_z^*$ act on $\p^1$ with weight 1 at $0\in \p^1$ and let 
\[ z:= e_{\BC_z^*}(\BC_{\mathrm{std}}) \]
be the class of the weight 1 representation $\BC_{\mathrm{std}}$ in the equivariant cohomology of a point. 
 The space  $V(M(\vsf),\beta)$ inherits $\BC_z^*$-action. Since $V(M(\vsf),\beta)$ is not proper, we define the virtual fundamental class of $V(M(\vsf),\beta)$ by the torus localisation, 
\[[V(M(\vsf),\beta)]^{\mathrm{vir}}:=\frac{[V(M(\vsf),\beta)^{\BC_z^*}]^{\mathrm{vir}} }{e_{\BC_z^*}(\CN^{\mathrm{vir}})}\in H_{*}(V(M(\vsf),\beta)^{\BC_z^*})[z^\pm].\]
We are ready to define (Donaldson-Thomas) \textit{Vertex function}, also known as \textit{I-function} in the GIT quasimap theory . 
\begin{defn}[Vertex function and truncated Vertex function] \
	\begin{align*}
		\mathsf{V}_\beta(z)&:= \evsf_* [V(M(\vsf),\beta)]^{\mathrm{vir}} \in H^*(M(\vsf))[z^\pm] \\
	 \mu_{\beta}(z)&:=[z\mathsf{V}_\beta]_{z^{\geq 0}}  \in H^*(M(\vsf))[z].
	\end{align*}
	\end{defn}
If $M(\vsf)=(\BC^2)^{[d]}$, the function $\mathsf{V}_\beta(z)$ is the classical 1-leg Vertex function from \cite{MNOP1, MNOP}, hence the origin of its name. Vertex is of fundamental importance for many reasons, among which are the following ones:
\begin{itemize}
\item  it is a building block for Donaldson-Thomas theory, \cite{MOOP},\cite{AKMV};
\item  it solves quantum differential equations and the quantum Knizhnik–Zamolodchikov equations, \cite{MO};
 \item it has an expression in terms of flags of sheaves, thereby connecting us to yet  another enumerative theory, \cite{Ob};
 \end{itemize}

For us, it will serve a different purpose - it will be responsible for the quasimap wall-crossing. 

\begin{thm}[Quasimap wall-crossing,  \cite{N}]\label{thetheorem} If $(g,n)\neq (0,0) ,(0,1)$, then 
\begin{multline*}
\langle \lambda_1 \psi^{k_1}_1, \dots,\lambda_n \psi^{k_n}_n \rangle^{-}_{g,\beta} = \langle \lambda_n \psi^{k_1}_1, \dots,\lambda_1 \psi^{k_n}_n \rangle^{+}_{g,\beta} \\
+\sum_{\underline{\beta}} \langle \lambda_n \psi^{k_1}_1, \dots,\lambda_n \psi^{k_n}_n, \mu_{\beta_1}(-\psi_{n+1}), \dots, \mu_{\beta_k}(-\psi_{n+k}) \rangle^{+}_{g,\beta_0}/k!,
\end{multline*}
where $\underline{\beta}=(\beta_0, \beta_1, \dots, \beta_k)$ such that $\beta=\sum^{i=k}_{i=0}\beta_i$ and $\beta_i\neq 0$ for $i\geq 1$. 
	\end{thm}
\textit{Sketch of Proof.} We use Zhou's master space, \cite{YZ}. The space of $\epsilon$-stabilities, which is just $\BR_{>0}$, has a chamber-wall structure. There are finitely many walls,  and for each wall $\epsilon_0 \in \BR_{>0}$ one can construct a master space,
\[MQ^{\epsilon_0}(M(\vsf),\beta),\]
with a natural $\BC_z^*$-action. Let $\epsilon_+$ and $\epsilon_-$ be values close to $\epsilon_0$ from the right and from left respectively.   The $\BC_z^*$-fixed locus of the master space is roughly of the following form 
\begin{multline*}
 MQ_{g,n}^{\epsilon_0}(M(\vsf),\beta)^{\BC_z^*} \approx  Q_{g,n}^{\epsilon_-}(M(\vsf),\beta) \cup Q_{g,n}^{\epsilon_+ }(M(\vsf),\beta) \\ 
 \cup \coprod_{\underline{\beta}} Q_{n+k}^{\epsilon_+}(M(\vsf),\beta_0)\times_{M(\vsf)^k} \prod_{i=1}^{i=k}V(M(\vsf), \beta_i) ,
 \end{multline*}
where the union is taken over $\underline{\beta}$, such that $\deg(\beta_i)=1/\epsilon_0$ for $i\geq 1$. In reality, the fixed locus has this expression up to certain modifications, which have a minor effect on the enumerative geometry. After taking the residue  in the localisation formula, we obtain a relation between the associated classes. This relation gives the wall-crossing formula after crossing all the walls on the way from $\epsilon=-$ to $\epsilon=+$.
\qed
\\

The full Vertex function $\mathsf{V}_\beta(z)$ also plays a role in genus 0 wall-crossing, we refer to \cite[Theorem 6.6]{N} for a precise statement. 

	

\subsection{Applications of the wall-crossing} \label{applications1}
In what follows, we give more freedom to $M(\vsf)$,  allowing it also to be either of the following:
\begin{itemize}
	\item a  moduli space of sheaves\footnote{For example, defined via moduli spaces of framed sheaves on some compactification.} on a quasi-projective surface with an action of a torus $T$, whose fixed locus is proper, e.g.\ Hilbert scheme of points on $\BC^2$; 
	\item a moduli space of objects in a tilted heart, e.g.\ $S^{[d]} \subset \CS_\#^{[d]}$. 
	\end{itemize}
The results extend to these cases, as long as three properties from the proof of Theorem \ref{main} can be checked (in all of the examples below it is true by arguments from \cite{N}).  

\subsubsection{Symplectic surfaces in general}
If $S$ carries a holomorphic symplectic form, we either require that the torus action scales the symplectic form or that the obstruction theory can be reduced for all $\epsilon$-stable quasimaps. 
\begin{cor} \label{important}If $S$ is holomorphic symplectic\footnote{One has to be careful here, because if $S$ is not projective, we have to consider framed sheaves, then there are some situation-dependent technicalities \cite{Sala} regarding the symmetry of the obstruction theory. An easier approach would be to forget about quasimaps and consider just threefolds $S\times C$. In this case, we essentially need our objects to be compactly supported (like stable pairs) and $S$ to have a trivial canonical bundle.} and $\evsf \colon V(M(\vsf), \beta) \rightarrow M(\vsf)$ is proper, then 
	\[\langle \lambda_1 , \dots,\lambda_n  \rangle^{-}_{g,\beta} = \langle \lambda_1, \dots,\lambda_n \rangle^{+}_{g,\beta},\] 
	unless $(g,n)=(0,1), (0,2), (1,0),(0,0)$. 
	\end{cor}
\textit{Proof.} In this case, the obstruction theory of $V(M(\vsf), \beta)$ has a cosection. The (reduced) virtual dimension of $V(M(\vsf), \beta)$ can be readily calculated, 
\[\text{red. dim}(V(M(\vsf), \beta))= \dim(M(\vsf))+1.\]
Since $\evsf$ is proper, the Vertex function
\[ \mathsf{V}_\beta(z) \in H_T^*(M(\vsf))[z^\pm]   \]
does not have poles in $T$-weights (i.e.\ we do not need to localise to pushforward the class). Hence by the dimension constraint, it has the following form (after taking reduction),
\[ \mathsf{V}_\beta(z)= \mu_\beta \mathbb 1/z+O(1/z^2), \]
where $\mu_\beta \in \BQ$ and $\mathbb 1$ is the fundamental class of $M(\vsf)$. Using  Theorem \ref{thetheorem} and the string equation, we obtain the claim. 
\qed
\begin{example} We start with a non-example. The pair $(\BC^2)^{[d]} \subset (\CC^2)^{[d]}$  does not satisfy the requirements of Corollary \ref{important}, because $\evsf$  is not proper. Indeed,  points of type (2) can "escape to infinity" (cf. Section \ref{points}).  In particular, the Vertex function $\mathsf{V}_\beta(z)$ has higher powers of $z$, which makes Theorem \ref{thetheorem} harder to apply. 
	\end{example}
\begin{example}  Consider now the perverse pair $(\BC^2)^{[d]} \subset (\CC^2)^{[d]}_\#$, then $\evsf$ is proper, because now points of the type (2) cannot "escape to infinity" and there are no points of type (3). In particular, we obtain that 
	\[ \mathsf{GW}_{0,3}((\BC^2)^{[d]})= \mathsf{PT}(\BC^2\times \p^1/ \BC^2_{0,1,\infty} ),\] 
	where $\BC^2_{0,1,\infty}=\BC^2\times\{0,1,\infty\}$. This is a result from \cite[Theorem 6]{OP10}. 
	
	More generally, consider an ADE surface $S$. The evaluation map $\evsf$ for the pair $(S^{[d]} ,\CS_{\#}^{[d]})$  is also proper, since all proper curves in $S$ are rigid. Hence we obtain 
	\[ \mathsf{GW}_{0,3}(S^{[d]})= \mathsf{PT}(S\times \p^1/ S_{0,1,\infty}).\] 
	\end{example}
\begin{example} For an ADE surface $S$, consider now the following torsion pair 
	\begin{align*}
		\CT&=\{A \in \Cohc(S) \mid \dim(A)\leq 1\}  \\
		\CT^{\perp}&=\{B \in \Cohc(S) \mid \Hom(A,B)=0 \ \forall A \in \CT\},
	\end{align*}
and the corresponding titled heart $\Cohc(S)_{\star}$ with the  associated moduli space of objects $\CS_{ \star}^{[d]}$ in the class $(1,0,-d)$. The evaluation map $\evsf$ is also proper for the pair $S^{[d]} \subset \CS_{\star}^{[d]}$. In this case, moduli spaces of quasimaps are isomorphic to moduli spaces of Bryan--Steinberg pairs \cite{BS}. Hence we obtain 
	\[  \mathsf{PT}(S\times \p^1/ S_{ 0,1,\infty})=\mathsf{GW}_{0,3}(S^{[d]})= \mathsf{BS}(S\times \p^1/ S_{0,1,\infty}),\] 
a result similar to \cite{Liu}. 
	\end{example}
\vspace{0.4cm}
\noindent\textbf{Observation.} Examples above suggest  that properness of the evaluation map $\evsf$ can be achieved for any moduli space $M(\vsf)$ after an appropriate tilt.
\\

\begin{example} Moduli spaces of sheaves on a $K3$ surface $S$ is ideal for Corollary \ref{important}, because the evaluation map $\evsf$ is always proper. As in the examples above, we get
\[   \mathsf{DT}(S\times \p^1/ S_{0,1,\infty})=\mathsf{GW}_{0,3}(S^{[d]})= \mathsf{PT}(S\times \p^1/ S_{0,1,\infty}).\] 
Moduli spaces of stable higher-rank sheaves on a $K3$ surface are deformation equivalent to Hilbert schemes (together with a chosen algebraic curve class). Hence by deformation invariance of Gromov--Witten invariants, we obtain 
\[ \mathsf{DT}_{\rk>1}(S\times \p^1/ S_{0,1,\infty})\cong\mathsf{GW}_{0,3}(S^{[d]})=\mathsf{DT}(S\times \p^1/ S_{0,1,\infty}), \] 
a result explained in more detail in \cite{NK3}. 
\end{example}
\subsubsection{More on K3 surfaces}
Corollary \ref{important} is powerful, however it does not apply to very  important invariants - genus-0 2-point invariants, which includes genus-1 0-point invariants as a special case. In this case, we need to compute the truncation of  Vertex  $\mu_\beta(z)$. For $\BC^2$, the entire 1-leg Vertex was computed in \cite{MNOP1, MNOP}, see also \cite{OP10, MOOP}. For a $K3$ surface, $\mu_\beta(z)$ was computed in \cite{Ob}. As a corollary of this computation and Theorem \ref{thetheorem}, we get the wall-crossing part of the Igusa cusp form conjecture \cite{OPa}.
\begin{cor} If $S$ is a K3 surface, then
	\[ \mathsf{GW}_E(S^{[d]})+F(p,q,\tilde q)=\mathsf{DT}(S\times E)=-\frac{1}{\chi(p,q,\tilde q)},\]
	where $\chi(p,q,\tilde q)$ is the Igusa cusp form and $F(p,q,\tilde q)$ is detetmined in \cite{Ob}. The second equality is determined in \cite{OPi}. 
\end{cor}
In the same vein, the wall-crossing was one of the ingredients for the proof of Holomorphic anomaly equation for $K3$ geometries in \cite{Ob22}.

\subsubsection{Del Pezzo surfaces} \label{Delpezzo}
Finally, if $S$ is a del Pezzo surface, then the truncation $\mu_\beta(z)$ was computed for the pair $S^{[n]}\subset \CS^{[n]}_\#$ in \cite{N}. To state the result, let us use the Nakajima basis (assuming $d>1$), 
\[(\gamma,\mathsf{k})\in H_2(S) \oplus \BZ E =H_2(S^{[d]},\BZ)\]
where $E$ is the primitive class of a curve contracted by the Hilbert-Chow map. Let us also define the following divisor class, 
\[\mathrm{c}_1(S)_d:=\mathsf{p}_{-1}(c_1(S)) \cdot \mathsf{p}_{-1}^{d-1}(\mathbb 1) \cdot \mathbb 1_S \in H^2(S^{[d]},\BZ).\]
The truncated Vertex function then takes the following form, 

\[ \sum_{(\gamma,\mathsf{k})} \mu_\beta(z)x^\gamma y^\mathsf{k}=\log (1+y)\mathrm{c}_1(S)_d.\]
Using Theorem \ref{thetheorem},  the divisor equation and the following intersection number for a class $(\gamma,\mathsf{k})$, 
\[ (\gamma,\mathsf{k})\cdot\mathrm{c}_1(S)_d=\gamma  \cdot \mathrm{c}_1(S),\]
 we obtain the following result. 
\begin{cor}\ \label{important2}If $S$ is a del Pezzo surface, then 
		\[ \sum_{\msf\geq 0} \langle \lambda_1 , \dots,\lambda_n  \rangle^{-}_{g,(\gamma,\mathsf{k})} y^{\mathsf{k}} = (1+y)^{\gamma\cdot\mathrm{c}_1(S)}\cdot \sum_{\msf\geq 0}  \langle \lambda_1, \dots,\lambda_n \rangle^{+}_{g,(\gamma,\mathsf{k})}y^{\mathsf{k}},\] 
		unless $(g,n)=(0,1), (0,2), (1,0),(0,0)$. 
	\end{cor}
\vspace{0.1cm}
\section{Enumerative mirror symmetry for moduli spaces of Higgs bundles}
We will now apply the theory of quasimaps to moduli spaces of Higgs bundles on a curve. Quasimaps give an enumerative realisation of Kapustin--Witten's considerations, \cite{KW}. In particular, by using the correspondence between  quasimaps and Vafa--Witten theory, we get a statement which relates curve counts in a moduli space of Higgs $\mathrm{SL}_\rsf$-bundles and curves counts in a moduli space of Higgs $\mathrm{PGL}_\rsf$-bundles. 

\subsection{Preliminaries } Let $\rsf,\dsf \in \BZ_{\geq0}$, such that $\rsf$ is prime and $0 \leq \dsf<\rsf$. Let $X$ be a smooth curve of genus $g\geq 2$. A Higgs sheaf $X$ is a pair $(F,\phi)$, a sheaf $F$ on $X$ and a morphism $\phi \in \Hom(F,F\otimes \omega_X)$. We define spaces of Higgs $\mathrm{SL}_\rsf$-sheaves, 
\begin{align*} \FMs&(\dsf)=\left\{ \text{Higgs sheaves }(F,\phi)   \text{ on } X \  \Bigl\rvert \ \arraycolsep=0.1pt\def\arraystretch{1} \begin{array}{c}  \det(F)=L, \tr(\phi)=0, \\[.001cm]  \ch(F)=(\rsf,\dsf) \end{array} \right\} \\
  &\rotatebox[origin=c]{90}{$\subset$}   \\ 
\Ms&(\dsf)=\{\text{slope semistable Higgs $\mathrm{SL}_\rsf$-bundles} \}.
\end{align*}

 The finite group of $\rsf$-torsion line bundles on $X$, denoted by $\Gamma_X$, acts on $\FMs(\dsf)$ by tensoring sheaves with lines bundles. We define the dual $\mathrm{PGL}_\rsf$-spaces, 
\[\FMp(\dsf):= [\FMs(\dsf)/\Gamma_X] \quad \text{and} \quad \Mp(\dsf):=[\Ms(\dsf)/\Gamma_X].\] 

We will need two classes, the class of the theta line bundle on $\FMs(\dsf)$ (the ample generator of the rational Picard group) and the class of the  $\BZ_\rsf$-gerbe of  SL-liftings of the universal family on $\FMs(\dsf)$ (the class used by \cite{HT}), denoted as 
\[ \Thetas \in H^2(\FMs(\dsf),\BZ) \quad \text{and} \quad \alphas \in H^2(\FMs(\dsf), \BZ_\rsf), \]
respectively. Both classes descend to classes on $\FMp(\dsf)$, 
\[ \Thetap \in H^2(\FMp(\dsf),\BZ)\quad \text{and} \quad \alphap \in H^2(\FMp(\dsf), \BZ_\rsf). \]
Finally, by scaling the Higgs field, we have torus actions on both spaces, 
\[ \BC^*_t \curvearrowright  \FMs(\dsf)\quad \text{and} \quad \BC^*_t  \curvearrowright  \FMp(\dsf) . \] 
Now that we have set up the stage, we will write $\Msp(\dsf)$ to mean either $\Mp(\dsf)$ or $\Mp(\dsf)$, the same convention applies to any other notation.
\subsection{Quasimaps, $\dsf\neq0$} Assume for the moment that $\dsf\neq 0$, i.e.\ there are no strictly semistable sheaves.  Let $E$ be an elliptic curve. 

\begin{defn} We say that a quasimap $f \colon E \dashrightarrow \Msp(\dsf)$ is of degree $(\wsf,\asf) \in \BZ \oplus \BZ_\rsf$, if 
	\vspace{-0.2cm}
\begin{align*}f^*\Thetasp&=\wsf \in H^2(E,\BZ)\cong \BZ, \\
	 f^*\alphasp&=\asf \in H^2(E, \BZ_\rsf)\cong \BZ_\rsf. 
	 \end{align*}
	 \end{defn}
 For $\wsf\neq0$,  let  
\[ Q_{E}(\Msp(\dsf), \asf,\wsf)^{\bullet} \] 
be the moduli spaces of quasimaps from $E$ to $\Msp(\dsf)$ up to translations\footnote{We identify quasimaps, if they equal up to a translation of $E$.} of $E$. 
The (reduced) expected dimension of these spaces is 0, so we can define a virtual number of quasimaps from $E$, 

\begin{align*}
	 \QMsp^{\asf, \bullet}_{\dsf,\wsf}:= \int_{[Q_{E}(\Msp(\dsf),\asf,\wsf)^{\bullet}]^{\mathrm{vir}}} 1 \in \BQ,
\end{align*} 
where we secretly use virtual localisation with respect to $\BC^*_t$-actions.  
\subsection{Vafa-Witten theory} \label{vwqm} Applying a version of Theorem \ref{ident} to this situation, we obtain 
\begin{align*} 
	\begin{split}
  Q_{E}(\Ms(\dsf), \asf,\wsf)&\cong \{ \text{Higgs sheaves on } X\times E\} \\
Q_{E}(\Mp(\dsf), \asf,\wsf)&\cong \{ \text{Higgs-Azumaya algebras on } X\times E\}.
\end{split}
\end{align*}
Higgs sheaves on a surface is a subject of Vafa--Witten theory, mathematically defined in \cite{TT, TT2}. 
Taking identification by translations of $E$ corresponds to taking an insertion ($\mu$-insertion). 
 \subsection{Quasimaps, $\dsf=0$} \label{d0} If $\dsf=0$, then there are strictly semistable Higgs sheaves, therefore $M(0)$ is an Artin stack. Counting curves in an Artin stack requires special treatment, we invite the reader to contemplate on what it means to count curves in $B\mathrm{GL}_\rsf$.  In this work, we take the following route, which is in accordance with \cite[Section 6]{Wit}. We use the relation between quasimaps and sheaves from Theorem \ref{ident}, which allows us to impose an extra stability condition on quasimaps  $f\colon E \dashrightarrow M(0)$. 
 
 \begin{defn}
 A quasimap $f\colon E \dashrightarrow M(0)$ is \textit{semistable}, if the associated Higgs sheaf on $X\times E$ is semistable with respect to a polarisation of the form $\CO_X(1)\boxtimes \CO_E(k)$ for $k\gg 1$. 

\end{defn}
By \cite[Proposition 5.4]{Nhiggs}, the condition above is automatic for $\dsf\neq 0$; it is also independent of $k$, as long as $k\gg 1$. A similar definition can be given in the $\mathrm{PGL}_\rsf$  case. If $\asf\neq 0$, then all quasimaps are stable, the resulting moduli space is a scheme and we can define invariants in the usual way. If $\asf=0$,  we use the theory of Joyce--Song \cite{JS} and Joyce \cite{J} to define the quasimap invariants for $\Msp(0)$, 
\[\QMsp^{\asf, \bullet}_{0,\wsf} \in \BQ,\] 
we expect the virtual fundamental class of \cite{ArP} can be used instead. 

\subsection{Instanton vs.\  Monopole}
Consider  the $\BC_t^*$-action on $\Msp(\dsf)$ . We have a decomposition of the $\BC_t^*$-fixed loci, 
\begin{align*}
	\Msp(\dsf)^{\BC_t^*}&=\Nsp(\dsf) \cup \Msp(\dsf)^{\mathrm{nil}}, 
\end{align*}
where $\Nsp(\dsf)$ are moduli spaces of semistable bundles;  $\Msp(\dsf)^{\mathrm{nil}}$ are loci of semistable Higgs bundles with nilpotent Higgs fields. For every $\dsf$, the quasimap invariants can be split into instanton and monopole contributions according to this decomposition, 
\[\QMsp^{\asf, \bullet}_{\dsf,\wsf}=  \textsf{instanton}+ \textsf{monopole}, \]
where \textsf{instanton} and  \textsf{monopole} summands are the invariants associated to $\Nsp(\dsf)$ and  $\Msp(\dsf)^{\mathrm{nil}}$ respectively. Terminology is borrowed from Vafa--Witten theory, e.g.\ see \cite{TT, TT2}. 
\subsection{Enumerative mirror symmetry} \label{applications2} Using identification from Section \ref{vwqm}, we can translate physical calculations of  \cite{MM} into mathematical conjectures which completely determine the genus 1 theory of $\FMsp(\dsf)$. We state conjectures for $\rsf=2$, see \cite[Section 7.2]{N} for an arbitrary prime rank.

We put invariants together into generating series,
\begin{align*}
	\QMsp^{\asf}_{\dsf}(q):= \sum_{\wsf>0} \QMsp^{\asf, \bullet}_{\dsf, \mathsf w}q^{\mathsf w}.
\end{align*}
Let $\tilde{\eta}(q)=\prod^{\infty}_{k=1}(1-q^k).$ We then define
\begin{align*}
\textsf{monopole}  & \begin{cases}	&\hspace{0.05cm}\Usf_1(q)= \log \tilde{\eta}(q^4)\end{cases}\\
	\textsf{instanton} &\begin{cases}  &\Usf_2(q)= \log \tilde{\eta}(q)\\
&\Usf_3(q)= \log\tilde{\eta}(-q). \end{cases}
\end{align*}
The invariants have the following form. 
\begin{conjno}[\hspace{-0.1cm} \cite{Nhiggs}] \label{wnonzero} If $g\geq 2$ and $\rsf=2$, then
	\begin{align*}
		\QMs^0_{\dsf}(q)&=(-1)^\dsf(2-2g)2^{4g-1}  \Usf_1(q) \\	
		\QMs^1_{\dsf}(q)&=(2-2g) 2^{2g-1 }(\Usf_2(z)+(-1)^\dsf\Usf_3(q)) \\
		\\
		\QMp^0_{\dsf}(q)&=(-1)^{\dsf}(2-2g)2^{2g-1}\Usf_1(q) \\
		\QMp^1_{\dsf}(q)&=(2-2g) (2^{4g-1}  \Usf_2(q) +(-1)^{\dsf} 2^{2g-1}  \Usf_3(q)).
	\end{align*}
\end{conjno}
There is a simple symmetry between these invariants, which we call \textit{Enumerative mirror symmetry}.  It partially interchanges monopole and instanton contributions.  To state it, let the symmetric group $S_3$ act on the linear span $ \BQ\langle \Usf_i(q) \rangle $ by permutation, 
\[ \sigma \cdot \Usf_i(q):=\Usf_{\sigma(i)}(q), \quad \sigma \in S_3.\]
This permutation is essentially induced by the transformation $\tau \rightarrow -1/\tau$ for $q=e^{\pi i \tau}$. 
\begin{cor}[Enumerative mirror symmetry] If $g\geq 2$ and $\rsf=2$, then
	
	\[(12)\cdot 2\QMs^{\asf}_{\dsf}(q)=\sum_{\dsf'} \sum_{\mathsf a'} (-1)^{  \mathsf d \cdot \mathsf a'+\dsf'\cdot \asf' } \hat{\mathsf{QM}}^{\asf'}_{\dsf'}(q).\]
	
\end{cor}
In \cite{Nqm}, we prove a part of the Conjecture, evoking Theorem \ref{thetheorem} and a curios relation between quasimaps from $E$ to $\Ms(\dsf)$ and quasimaps from $X$ to a moduli spaces of Higgs bundles on $E$. 
\begin{thm}[ \hspace{-0.1cm}\cite{Nqm}] If $g\geq 2$ and $\rsf=2$, then 
\[	\QMs^1_{\dsf}(q)=(2-2g) 2^{2g-1 }(\Usf_2(q)+(-1)^\dsf\Usf_3(q)). \]
	\end{thm}
Moreover, by Theorem \ref{thetheorem} and the analysis of the truncated Vertex functions from \cite[Section 10]{Nhiggs}, we obtain the following. 
\begin{thm} If $\dsf\neq 0$ and $\gcd(\rsf,\wsf)=1$, then 
	\[\QMs^{\asf,\bullet}_{\dsf,\wsf}=\GWs^{\asf,\bullet}_{\dsf,\wsf},\]
	where $\GWs^{\asf,\bullet}_{\dsf,\wsf}$ are Gromov-Witten genus 1 invariants. 
	\end{thm}

\section{Gromov-Witten/Hurwitz wall-crossing} 

 Given a smooth projective variety $X$. We want to relate Gromov--Witten theory of the orbifold symmetry product $X^{(d)}=[X^d/S_d]$ to Gromov-Witten theory of $X\times C$ in the same way it was done with quasimaps. Compare  Figure \ref{GWH} to Figure \ref{summary}.
 
	\begin{figure} [h!]	\vspace{1.5cm}
	\scriptsize
	\[
	\begin{picture}(200,0)(-30,-30)
		\thicklines
		\put(26,0){\line(1,0){29}}
		\put(95,0){\line(1,0){25}}
		\put(75,0){\makebox(0,0){\textsf{GW Vertex}}}
		\put(143,0){\makebox(0,0){ $\mathsf{GW}_C(X^{(d)})$}}
		\put(0,0){\makebox(0,0){$\mathsf{GW_{rel}}(S\times C)$}}
	\end{picture}
	\vspace{-0.5cm}	\]
	\caption{Gromov--Witten/Hurwitz wall-crossing}
	\label{GWH}

\end{figure}

\subsection{Unramfied maps} Let us firstly translate Gromov--Witten theory of $X^{(d)}$ into a language more appropriate for our purposes. The data\footnote{More precisely, $X^{(d)}$ is an orbifold, so we have to consider twisted curves, however, we will not spend time on it, as we will use a  different point of view anyway. } of a map $f \colon (C,p_i) \rightarrow X^{(d)}$ is equivalent to data of the following correspondence

\[
\begin{tikzcd}[row sep=small, column sep = small]
	P \arrow[r, "f_{X}"] \arrow{d}[swap]{f_{C}} & X  & \\
	(C,p_i) & & 
\end{tikzcd}
\]
such that $f_C$ is an unramified map (local isomorphism) of degree $d$. More precisely, it is unramified everywhere except at marked points and nodes:
\begin{itemize} \item at a marked point, a ramification profile is fixed,
	\item  at a node $\BC[x_1,x_2]/(x_1x_2)$, it is of the form $(x_1,x_2) \mapsto (x^k_1,x_2^k)$.
	\end{itemize} 
The two  conditions are called  \textit{admissibility}, all maps in this section are assumed to satisfy it. Unramified maps from a curve to another curve is a subject of Hurwitz theory, hence the following analogy 
\[ \textit{Gromov--Witten theory of }X^{(d)} = \textit{Hurwitz theory with a target }X. \]

\begin{figure}[h!]
	\centering
	\includegraphics[scale=0.25]{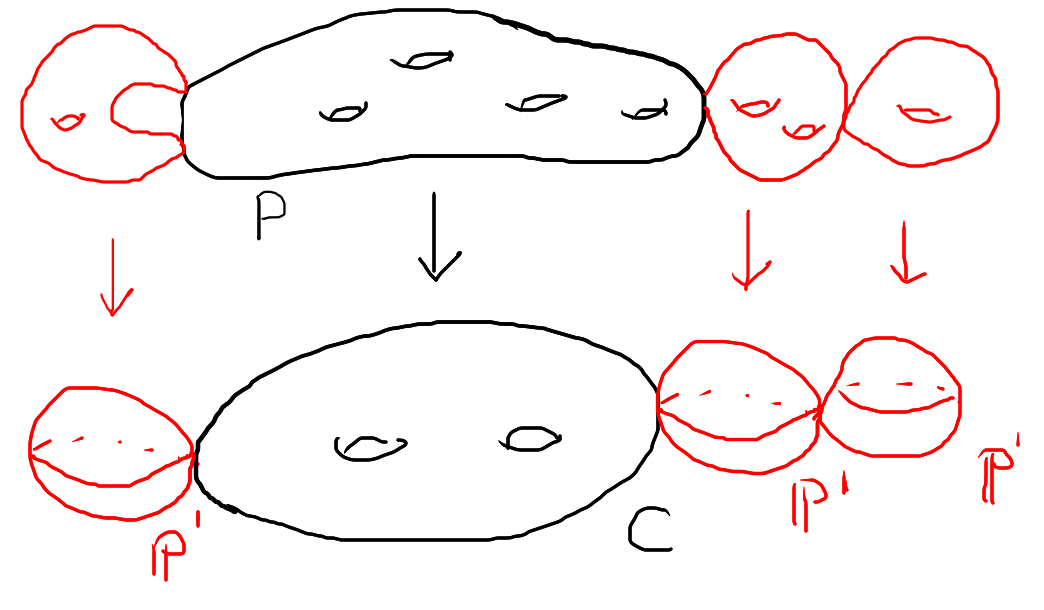}
	\caption{Unramified maps}
	\label{unramified maps}
\end{figure}

 By Riemann--Hurwitz formula, the  arithmetic genus $g(P)$ of the source curve is determined by $d$, $g(C)$ and the ramification profile of marked points. However, we want  $g(P)$ to have some independence from these discrete variables. To achieve it, we allow $f_C$ to have simple branching,  i.e.\ points $b_j \in C$, such that $| f_C^{-1}(b_j)|=d-1$. We treat the number of $b_j$ as an additional discrete degree variable, it is related to $g(P)$ by  Riemann--Hurwitz formula (as long as other discrete variables are fixed), 
\[ |b_j| \iff g(P).\]
We put the degree of $f_X$ and $|b_j|$ together into a \textit{refined degree}, denoting them by $\beta=(\gamma,\msf) \in H_2(S)\oplus \BZ$. 
\pagebreak
 Summing up the discussion above, we introduce moduli spaces of pairs $(f_C,f_X)$, such that $f_C$ is unramified,  
\[ \Mbar_{g,n}(X^{(d)},\beta)=\left\{\text{unramified }(f_C,f_X) \  \Bigl\rvert \ \arraycolsep=0.1pt\def\arraystretch{1} \begin{array}{c}  \deg(f_X)=\gamma, |b_j |=\msf,  \\[.001cm] g(C)=g , | p_i |=n  \end{array}  \right\}.\] 
This moduli space can essentially be understood\footnote{See \cite[Section 2.3]{N22} for the precise statement.} as a moduli space of maps to $X^{(d)}$ with a refined degree. We also suggest the reader to treat Figure \ref{unramified maps} and Figure \ref{bubbles} as equal, and please forgive the author, if Figure \ref{unramified maps} contradicts Riemann--Hurwitz formula. 
\vspace{0.3cm}
\vspace{-0.3cm}

	\subsection{Ramified maps }	
The other extreme is when $f_C$ is allowed to have arbitrary ramifications and contracted components. 
\begin{defn} Abusing terminology, we call a map $f_C\colon P \rightarrow C$ \textit{ramified}, if it has one (or more) of the following:
	\begin{itemize}
		\item ramification of order greater than 2; 
		\item contracted components; 
		\item nodes mapping to regular locus,
		\end{itemize}
	\end{defn}

Hence in this case, we just have stable maps from $P$ to $X\times C$ of degree $(d,\gamma)$, as there are no restrictions on $f_C$ (apart from being admissible at marked points and nodes). In this case, we can prohibit $C$ to have rational tails, like we did with quasimaps in Definition \ref{stableq}. Hence we obtain a moduli a space of \textit{ramified} maps to $X^{(d)}$, 

\[ r\Mbar_{g,n}(X^{(d)}, \beta):= \Mbar_\msf(X\times \CC_{g,n}/ \CM_{g,n},(\gamma ,d )),\]
where  $\CC_{g,n}  \rightarrow \CM_{g,n}$ is the universal curve of the moduli space of semistable curves (no rationals tails). The integer $\msf$ is degree of the branching divisor,  we refer to the next section for more details on the branching divisor.

\begin{figure}[h!]\vspace{0.3cm}
	\centering
	\includegraphics[scale=0.25]{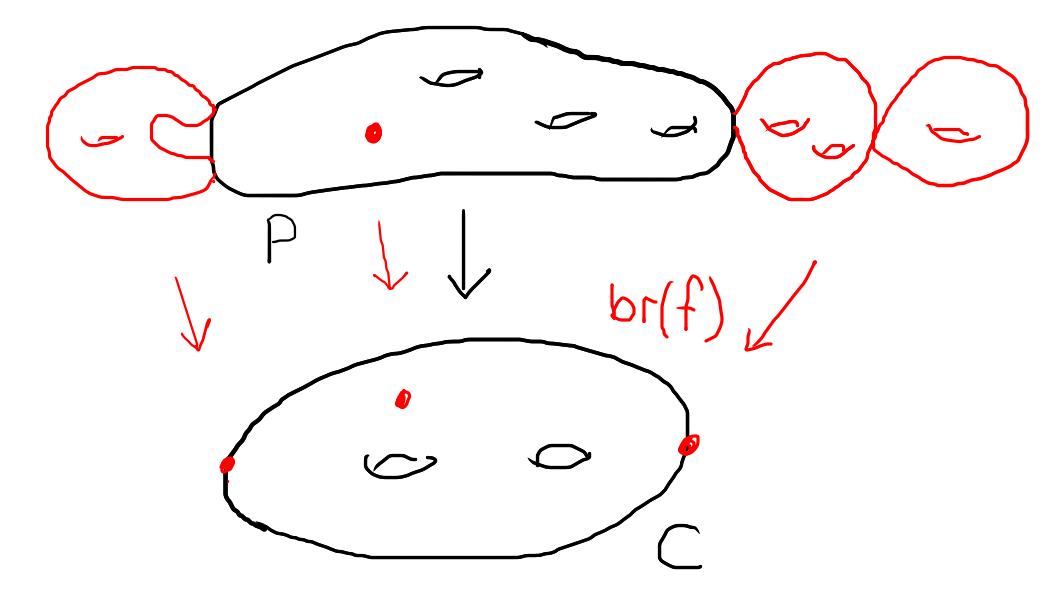}
	\caption{Ramified maps}
	\label{ramified maps}
	\vspace{-0.3cm}
\end{figure}

This moduli space is analogous to the moduli space of quasimaps for Hilbert scheme of $d$-points $S^{[d]}$. In fact, this analogy is very much justified, because a ramified map defines a rational map $f \colon C \dashrightarrow X^{(d)}$. The difference is that there is no ambient stack to keep track of the base points, hence contracted components  and ramifications replace them. Again, we suggest the reader to treat Figure \ref{ramified maps} and Figure \ref{Degenerate} as equal. 

\subsection{$\epsilon$-ramified maps}
There is a natural object that measures the difference between unramified and ramified maps. It is given by the branching complex ,
\[\br(f) =Rf_{C*}[f^*_C\Omega_C \rightarrow \Omega_P],\]
which is a skyscraper complex supported on branching points of $C$ (we can ignore what happens at markings and nodes). We will identify $\br(f)$ with the divisor given by the support of the complex weighted by its Euler characteristics. One can check that contracted components, ramification and nodes contribute to $\br(f)$ in the following way: 
\begin{itemize}
	\item ramifications of order $e$ contribute $e-1$;
	\item contracted components of genus $g$ contribute $2g-2$; 
	\item nodes mapping to the regular locus contribute $2$.  
	\end{itemize}
We can therefore interpolate between ramified and unramified, using the idea of $\epsilon$-\textit{ramification}. For that, we choose an ample line bundle $L$ on $X$ (of high enough degree). 
\vspace{0.2cm}

\noindent \textbf{Important.} Our terminology here slightly differs from the one of \cite{N22}, where we call such stability $\epsilon$-\textit{admissibility}. 
\begin{defn} Given $\epsilon \in \BR_{>0}$.  A map $(f_C,f_X)$ is $\epsilon$-ramified, if 
	\begin{itemize}
	\item $\forall$ branching points $b \in C, \ \mathrm{mult}_b\br(f)+\deg(f^*_XL_{|f_C^{-1}(b)})\leq e^{1/\epsilon}$,
	\item $\forall$ rational tails $R \subset C , \ \deg( \br(f)_{|R})+\deg(f^*_XL_{|f_C^{-1}(R)}) > e^{1/\epsilon}$,
	\item $|\Aut(f_C, f_X) | < \infty$. 
	\end{itemize}
	\end{defn} 
One can readily verify that 
\begin{align*}
+\text{-ramified maps } &= \text{ unramified maps } \\
- \text{-ramified maps } &= \text{ ramified maps}
\end{align*}
We define moduli spaces of $\epsilon$-ramified maps
\[ r\Mbar^\epsilon_{g,n}(X^{(d)},\beta)= \left\{\epsilon\text{-ramified }(f_C,f_X) \  \Bigl\rvert \ \arraycolsep=0.1pt\def\arraystretch{1} \begin{array}{c}  \deg(f_X)=\gamma, \deg(\br(f))=\msf,  \\[.001cm] g(C)=g , | p_i |=n  \end{array} \right\}\]
\begin{thm}[\hspace{-0.01cm}\cite{N22}] If $X$ is a smooth projective variety, then a moduli space $r\Mbar^\epsilon_{g,n}(X^{(d)}, \beta)$ is a proper Deligne-Mumford stack with a perfect obstruction theory. 
	\end{thm}

\subsection{$\epsilon$-ramified invariants} Recall the inertia stack of $X^{(d)}$,
\[ \CI X^{(d)}:= \coprod_{\eta \vdash d } [X^{\eta}/\Aut(\eta)], \]
where $\eta$ is a partition of $d$ of  length $l(\eta)$ and $X^{\eta}:=X^{l(\eta)}$.
There exist usual structures needed to define Gromov-Witten type invariants:
\begin{itemize}
	\item virtual fundamental class
	\[[r\Mbar^\epsilon_{g,n}(X^{(d)}, \beta)]^{\mathrm{vir}}\in H_{\mathrm{vdim}}(r\Mbar^\epsilon_{g,n}(X^{(d)}, \beta)); \]
	\item evaluation maps 
	\[\ev_i: r\Mbar^\epsilon_{g,n}(X^{(d)}, \beta) \rightarrow \CI X^{(d)}, \quad f \mapsto f_X(f_C^{-1}(p_i));\]
	\item $\psi$-classes, 
	\[ L_{i}=T^*_{p_i}C , \quad \psi_i:=\mathrm{c}_1(L_i). \]
\end{itemize}
\begin{defn} For classes $\lambda_i \in H^*_{\mathrm{orb}}(X^{(d)}):=H^*(\CI X^{(d)})$, we define $\epsilon$-ramified invariants 
	\[\langle \lambda_1 \psi^{k_1}_1, \dots,\lambda_n \psi^{k_n}_n \rangle^{\epsilon}_{g,\beta}= \int_{[r\Mbar^\epsilon_{g,n}(X^{(d)}, \beta)]^{\mathrm{vir}}} \prod^{n}_{i=1} \ev^*_i \lambda_i \cdot \psi^{k_i}.\]
\end{defn}

\subsection{Vertex and wall-crossing}
As in the case of quasimaps, we consider ramified maps for a fixed $\p^1$ without ramifications\footnote{In fact, some ramifications are allowed, but we treat it as ramifications of markings} at $\infty \in \p^1$, 
\[ V(X^{(d)},\beta) := \left\{ \text{ ramified } (f_{\p^1},f_X)  \  \Bigl\rvert \ \arraycolsep=0.1pt\def\arraystretch{1} \begin{array}{c}  \deg(f_X)=\gamma, \deg(\br(f))=\msf,  \\[.001cm]\text{ admissible at } \infty\in \p^1  \end{array} \right\}.\] 
By definition, we have an evaluation map, 
\[ \evsf \colon V(X^{(d)},\beta) \rightarrow \CI X^{(d)}, \quad f \mapsto f_X(f_{\p^1}^{-1}(\infty)) .\]
Using the $\BC^*_z$-action on $\p^1$, we define the virtual fundamental class by localisation, 
\[[V(X^{(d)},\beta)]^{\mathrm{vir}}:=\frac{[V(X^{(d)},\beta)^{\BC_z^*}]^{\mathrm{vir}} }{e_{\BC_z^*}(\CN^{\mathrm{vir}})}\in H_{*}(V(X^{(d)},\beta)^{\BC_z^*})[z^\pm].\] 
\begin{defn}[Gromov--Witten Vertex] \
	\begin{align*}
		\mathsf{V}_\beta(z)&:= \evsf_* [V(X^{(d)},\beta)]^{\mathrm{vir}} \in H^*_{\mathrm{orb}}(X^{(d)})[z^\pm]\\
		\mu_{\beta}(z)&:=[z\mathsf{V}_\beta]_{z^{\geq 0}}  \in H^*_{\mathrm{orb}}(X^{(d)})[z].
	\end{align*}
\end{defn}

\begin{thm}[Gromov--Witten/Hurwitz wall-crossing, \cite{N22}]\label{thetheorem2} If $(g,n)\neq (0,0) ,(0,1)$, then 
	\begin{multline*}
		\langle \lambda_1 \psi^{k_1}_1, \dots,\lambda_n \psi^{k_n}_n \rangle^{-}_{g,\beta} = \langle \lambda_n \psi^{k_1}_1, \dots,\lambda_1 \psi^{k_n}_n \rangle^{+}_{g,\beta} \\
		+\sum_{\underline{\beta}} \langle \lambda_n \psi^{k_1}_1, \dots,\lambda_n \psi^{k_n}_n, \mu_{\beta_1}(-\psi_{n+1}), \dots, \mu_{\beta_k}(-\psi_{n+k}) \rangle^{+}_{g,\beta_0}/k!,
	\end{multline*}
where $\underline{\beta}=(\beta_0, \beta_1, \dots, \beta_k)$ such that $\beta=\sum^{i=k}_{i=0}\beta_i$ and $\beta_i\neq 0$ for $i\geq 1$. 
\end{thm}

\subsection{Applications of the wall-crossings} \label{applications3} Firstly, combining Theorem \ref{thetheorem} and Theorem \ref{thetheorem2}, we obtain that Ruan's Crepant resolution conjecture \cite{Ru} and MNOP's relative Pandharipande--Thomas/Gromov--Witten correspondence \cite{MNOP,MNOP1} are equivalent up to a correspondence of truncated Vertex functions. 

\begin{cor} For a surface $S$, we have 
	\[ \mathsf{CRC} \ \mathrm{for} \ (S^{[d]}, S^{(d)})= \mathsf{relative \ PT/GW} \ \mathrm{for} \ S\times C \ + \ \mathsf{PT/GW}\  \mathrm{for} \ \mu_\beta(z).\]
	\end{cor}
Moreover, by \cite[Section 2.2]{BG} the change of variables $y=-e^{iu}$ is not necessary for the formulation of Crepant resolution conjecture - it arises due to failure of the divisor equation on $S^{(d)}$.

If $S$ is a del Pezzo surface,  then the truncated Vertex was computed in \cite[Proposition 5.1]{N22} and is equal to
\[ \sum_{(\gamma,\msf)}\mu_{(\gamma,\msf)}(z)x^\gamma u^\msf= \log\left(\frac{\sin(u/2)}{u/2}\right)\mathrm{c}_1(S)_d,\]
where the class $\mathrm{c}_1(S)_d$ corresponds to the one in Section \ref{Delpezzo}, we refer \cite[Section 5]{N22} for the precise definition of this class. Using Theorem \ref{thetheorem2} and the divisor equation\footnote{In general, divisor equations do not hold for orbifolds, but the class $\mathrm{c}_1(S)_d$ lies in the untwisted sector.}, we obtain the following result. 
\begin{cor}\ \label{important3}If $S$ is a del Pezzo surface, then 
	\[ \sum_{k\geq 0} \langle \lambda_1 , \dots,\lambda_n  \rangle^{-}_{g,(\gamma,\msf)} u^{\msf} = \left(\frac{\sin(u/2)}{u/2}\right)^{\gamma\cdot\mathrm{c}_1(S)}\cdot \sum_{k\geq 0} \langle \lambda_1, \dots,\lambda_n \rangle^{+}_{g,(\gamma,\msf)}u^{\msf},\] 
	unless $(g,n)=(0,1), (0,2), (1,0),(0,0)$. 
\end{cor}
Combining it with Corollary \ref{important2} and relative PT/GW from \cite{MOOP}, we prove some cases of Crepant resolution conjecture. 
\begin{cor} If $S$ is a toric del Pezzo surface, then  Crepant resolution conjecture holds for $(S^{[d]}, S^{(d)})$ for all curve classes in the case $(g,n)=(0,3)$.
	\end{cor}
	
	\bibliographystyle{amsalpha}
	\bibliography{notes}

\providecommand{\bysame}{\leavevmode\hbox to3em{\hrulefill}\thinspace}
\providecommand{\MR}{\relax\ifhmode\unskip\space\fi MR }
\providecommand{\MRhref}[2]{%
  \href{http://www.ams.org/mathscinet-getitem?mr=#1}{#2}
}
\providecommand{\href}[2]{#2}
\begin{thebibliography}{MNOP06b}

\bibitem[AKMV05]{AKMV}
M.~Aganagic, A.~Klemm, M.~Mari{\~n}o, and C.~Vafa, \emph{The topological
  vertex}, Commun. Math. Phys. \textbf{254} (2005), no.~2, 425--478.

\bibitem[AP22]{ArP}
D.~Aranha and P.~Pstragowski, \emph{{The Intrinsic Normal Cone For Artin
  Stacks}}, arXiv:1909.07478 (2022).

\bibitem[BG09]{BG}
J.~Bryan and T.~Graber, \emph{The crepant resolution conjecture}, Algebraic
  geometry, Seattle 2005. Proceedings of the 2005 Summer Research Institute,
  Seattle, WA, USA, July 25--August 12, 2005, Providence, RI: American
  Mathematical Society (AMS), 2009, pp.~23--42.

\bibitem[BS16]{BS}
J.~Bryan and D.~Steinberg, \emph{Curve counting invariants for crepant
  resolutions}, Trans. Am. Math. Soc. \textbf{368} (2016), no.~3, 1583--1619.

\bibitem[CKM14]{CFKM}
I.~{Ciocan-Fontanine}, B.~{Kim}, and D.~{Maulik}, \emph{{Stable quasimaps to
  GIT quotients}}, {J. Geom. Phys.} \textbf{75} (2014), 17--47.

\bibitem[Gro61]{Groth}
Alexandre Grothendieck, \emph{Techniques de construction et th{\'e}oremes
  d'existence en g{\'e}om{\'e}trie alg{\'e}brique. {IV}: {Les} schemas de
  {Hilbert}}, Sem. {Bourbaki} 13(1960/61), {No}. 221, 28 p. (1961)., 1961.

\bibitem[Gro85]{Gro}
M.~Gromov, \emph{Pseudo holomorphic curves in symplectic manifolds}, Invent.
  Math. \textbf{82} (1985), 307--347.

\bibitem[HT03]{HT}
T.~Hausel and M.~Thaddeus, \emph{Mirror symmetry, {Langlands} duality, and the
  {Hitchin} system}, Invent. Math. \textbf{153} (2003), no.~1, 197--229.

\bibitem[Joy21]{J}
D.~Joyce, \emph{{Enumerative invariants and wall-crossing formulae in abelian
  categories}}, arXiv:2111.04694 (2021).

\bibitem[JS12]{JS}
D.~Joyce and Y.~Song, \emph{A theory of generalized {Donaldson}-{Thomas}
  invariants}, Mem. Am. Math. Soc., vol. 1020, Providence, RI: American
  Mathematical Society (AMS), 2012.

\bibitem[Kon95]{Kon}
M.~Kontsevich, \emph{Enumeration of rational curves via torus actions}, The
  moduli space of curves. Proceedings of the conference held on Texel Island,
  Netherlands during the last week of April 1994, Basel: Birkh{\"a}user, 1995,
  pp.~335--368.

\bibitem[KW07]{KW}
A.~Kapustin and E.~Witten, \emph{Electric-magnetic duality and the geometric
  {Langlands} program}, Commun. Number Theory Phys. \textbf{1} (2007), no.~1,
  1--236.

\bibitem[Liu22]{Liu}
H.~Liu, \emph{{Equivariant K-theoretic enumerative invariants and wall-crossing
  formulae in abelian categories}}, arXiv:2207.13546 (2022).

\bibitem[Lur12]{Lur}
J.~Lurie, \emph{{Derived Algebraic Geometry XIV: Representability Theorems}},
  \url{https://www.math.ias.edu/~lurie/papers/DAG-XIV.pdf}.

\bibitem[MM21]{MM}
J.~Manschot and G.~W. Moore, \emph{{Topological correlators of $SU(2)$,
  $\mathcal{N}=2^*$ SYM on four-manifolds}}, arXiv:2104.06492 (2021).

\bibitem[MNOP06a]{MNOP1}
D.~Maulik, N.~Nekrasov, A.~Okounkov, and R.~Pandharipande,
  \emph{Gromov-{Witten} theory and {Donaldson}-{Thomas} theory. {I}}, Compos.
  Math. \textbf{142} (2006), no.~5, 1263--1285.

\bibitem[MNOP06b]{MNOP}
\bysame, \emph{Gromov-{W}itten theory and {D}onaldson-{T}homas theory. {II}},
  Compos. Math. \textbf{142} (2006), no.~5, 1286--1304. \MR{2264665}

\bibitem[MO19]{MO}
D.~Maulik and A.~Okounkov, \emph{Quantum groups and quantum cohomology},
  Ast{\'e}risque, vol. 408, Paris: Soci{\'e}t{\'e} Math{\'e}matique de France
  (SMF), 2019.

\bibitem[MOOP11]{MOOP}
D.~Maulik, A.~Oblomkov, A.~Okounkov, and R.~Pandharipande,
  \emph{Gromov-{Witten}/{Donaldson}-{Thomas} correspondence for toric 3-folds},
  Invent. Math. \textbf{186} (2011), no.~2, 435--479.

\bibitem[Nes21a]{N}
D.~Nesterov, \emph{Quasimaps to moduli spaces of sheaves}, arXiv:2111.11417
  (2021).

\bibitem[Nes21b]{NK3}
\bysame, \emph{{Quasimaps to moduli spaces of sheaves on a K3 surface}},
  arXiv:2111.11425 (2021).

\bibitem[Nes22]{N22}
\bysame, \emph{{Gromov--Witten/Hurwitz wall-crossing}}, arXiv:2208.00889
  (2022).

\bibitem[Nes23a]{Nhiggs}
\bysame, \emph{{Enumerative mirror symmetry for moduli spaces of Higgs bundles
  and S-duality}}, arXiv:2302.08379 (2023).

\bibitem[Nes23b]{Nqm}
\bysame, \emph{{On quasimap invariants of moduli spaces of Higgs bundles}},
  2023.

\bibitem[Obe21]{Ob}
G.~Oberdieck, \emph{{Multiple cover formulas for K3 geometries, wallcrossing,
  and Quot schemes}}, arXiv:2111.11239 (2021).

\bibitem[Obe22]{Ob22}
\bysame, \emph{{Holomorphic anomaly equations for the Hilbert scheme of points
  of a K3 surface}}, arXiv:2202.03361 (2022).

\bibitem[OP10]{OP10}
A.~Okounkov and R.~Pandharipande, \emph{The local {Donaldson}-{Thomas} theory
  of curves}, Geom. Topol. \textbf{14} (2010), no.~3, 1503--1567.

\bibitem[OP16]{OPa}
G.~Oberdieck and R.~Pandharipande, \emph{Curve counting on {$K3\times E$}, the
  {I}gusa cusp form {$\chi_{10}$}, and descendent integration}, K3 surfaces and
  their moduli, Progr. Math., vol. 315, Birkh\"{a}user/Springer, 2016,
  pp.~245--278.

\bibitem[OP18]{OPi}
G.~Oberdieck and A.~{Pixton}, \emph{{Holomorphic anomaly equations and the
  Igusa cusp form conjecture}}, {Invent. Math.} \textbf{213} (2018), no.~2,
  507--587.

\bibitem[PT09]{PT}
R.~{Pandharipande} and R.~P. {Thomas}, \emph{{Curve counting via stable pairs
  in the derived category}}, {Invent. Math.} \textbf{178} (2009), no.~2,
  407--447.

\bibitem[Rua03]{Ru}
Y.~Ruan, \emph{Discrete torsion and twisted orbifold cohomology}, J. Symplectic
  Geom. \textbf{2} (2003), no.~1, 1--24.

\bibitem[Sal12]{Sala}
F.~Sala, \emph{Symplectic structures on moduli spaces of framed sheaves on
  surfaces}, Cent. Eur. J. Math. \textbf{10} (2012), no.~4, 1455--1471.

\bibitem[TT18]{TT2}
Y.~Tanaka and R.~P. Thomas, \emph{Vafa-{Witten} invariants for projective
  surfaces. {II}: {Semistable} case}, Pure Appl. Math. Q. \textbf{13} (2018),
  no.~3, 517--562.

\bibitem[TT20]{TT}
\bysame, \emph{Vafa-{Witten} invariants for projective surfaces. {I}: {Stable}
  case}, J. Algebr. Geom. \textbf{29} (2020), no.~4, 603--668.

\bibitem[Wit10]{Wit}
E.~Witten, \emph{Mirror symmetry, {Hitchin}'s equations, and {Langlands}
  duality}, The many facets of geometry. A tribute to Nigel Hitchin, Oxford:
  Oxford University Press, 2010, pp.~113--128.

\bibitem[{Zho}22]{YZ}
Y.~{Zhou}, \emph{{Quasimap wall-crossing for GIT quotients}}, {Invent. Math.}
  \textbf{227} (2022), no.~2, 581--660.

\end{thebibliography}
\end{document}